\documentclass[11pt]{amsart}
\usepackage{amssymb,amsmath}
\usepackage{epsf,latexsym,graphicx}
\usepackage{psfrag}
\psfragscanon
\theoremstyle{plain}
\newtheorem{thm}{Theorem}
\newtheorem{prop}[thm]{Proposition}

\newtheorem{lemma}[thm]{Lemma}
\newtheorem{conjecture}[thm]{Conjecture}

\newtheorem{remark}[thm]{Remark}

\numberwithin{equation}{section}

\begin{document}

\title{Conjectures on the cohomology of the Grassmannian}
\author{Victor Reiner}
\author{Geanina Tudose}
\begin{abstract}
We give a series of successively weaker conjectures on the cohomology ring of the 
Grassmannian, starting with the Hilbert series of a certain natural filtration.  
\end{abstract}
\address{School of Mathematics, University of Minnesota, Minneapolis,
  MN, 55455 USA}
\address{The Fields Institute for Mathematical Sciences, 222 College Street, 
Toronto, ON, Canada}

\email{reiner@math.umn.edu}
\email{gtudose@fields.utoronto.ca}

\maketitle

\def\boxx{{\ell^k}}
\def\Gr{{\mathbb G}}
\def\CC{{\mathbb C}}
\def\QQ{{\mathbb Q}}
\def\ZZ{{\mathbb Z}}
\def\PP{{\mathcal P}}
\def\End{{\operatorname{End}}}
\def\Hilb{{\operatorname{Hilb}}}
\newcommand\qbinom[2]{ \left[ \begin{matrix} #1 \\ #2 \end{matrix} \right]_q }

\section{Introduction}

Let $\Gr(k,\CC^{k+\ell})$ denote the Grassmannian of $k$-planes in
$\CC^{k+\ell}$.  We pose a series of successively weaker conjectures on
its cohomology ring with rational coefficients. 
The first (Conjecture~\ref{filtration-hilb}) conjectures the 
Hilbert series for a natural filtration of the ring.  The last
(Conjecture~\ref{e-overkiller}) would greatly simplify the proof of
a result by Hoffman (Theorem~\ref{non-zero-alpha} below)
on the classification of its graded endomorphisms.

We denote this cohomology ring $R^{k,l}:=H^*(\Gr(k,\CC^{k+\ell}), \QQ).$
As a graded $\QQ$-algebra, $R^{k,l}$  has several natural descriptions
(see \cite{GH}, \cite{F}):
\begin{enumerate}
\item[(i)]
$$
R^{k,\ell} \cong \QQ[e_1,e_2,\ldots,e_k,h_1,h_2,\ldots,h_\ell]/
\left(\sum_{i+j=d} (-1)^i e_i h_j: d=1,2,\ldots,k+\ell\right).
$$
\item[(ii)]
$$
R^{k,\ell} \cong \QQ[e_1,e_2,\ldots,e_k]/(h_{\ell+1},h_{\ell+2},\ldots,h_{\ell+k})
$$
where here $h_r$ is interpreted as the {\it Jacobi-Trudi determinant}
$$
h_r : = \left|   
\begin{matrix}
e_1 &  e_2  & e_3  & \cdots &  & \\
1   &  e_1  & e_2  & \cdots &  & \\
0   &  1    & e_1  & \cdots &  & \\
\vdots & \ddots & \ & \vdots & & \\
0   &  \cdots    & 0    & \cdots & 1 & e_1
\end{matrix}
\right|
$$
\item[(iii)] $R^{k,\ell}$ has $\QQ$-basis given by the Schur
functions $s_\lambda$ for partitions $\lambda$ whose Ferrers diagram
fit inside the $k \times \ell$ box $\boxx$ (written $\lambda \subset \boxx$),
with multiplication given
by the usual {\it Littlewood-Richardson} structure constants
$$
s_\lambda s_\mu = \sum_{\nu \subset \boxx} c^\nu_{\lambda, \mu} s_\nu.
$$
Note the truncation to only those partitions $\nu$ satisfying $\nu \subset \boxx$.
\end{enumerate}
We give here a brief explanation of these various descriptions
and the connections between them.  The cohomology ring $R^{k,\ell}$ 
is generated by the Chern classes of the tautological $k$-vector bundle $\mathcal E$
on $\Gr(k,\CC^{k+\ell})$.  In description (i), $e_i$ represents $(-1)^i$ times the
$i^{th}$ Chern class $c_i(\mathcal E)$.  Here $h_j$ represents the $j^{th}$ Chern 
class $c_j({\mathcal F})$ for the quotient $\ell$-vector bundle 
${\mathcal F}$ that fits into the following 
exact sequence of bundles on $\Gr(k,\CC^{k+\ell})$:
$$
0 \rightarrow {\mathcal E} \rightarrow \CC^{k+\ell} 
\rightarrow {\mathcal F} \rightarrow 0
$$
in which the middle term a trivial $k+\ell$-vector bundle.
In description (ii), the $e_i$ are as in description (i).  One has
simply used the relations $\sum_{i+j=d} (-1)^i e_i h_j$ for
$d=1,\ldots,\ell$ to eliminate the redundant generators
$h_1,h_2,\ldots,h_\ell$, while
for $d=\ell+1,\ldots, \ell + k$,
these relations turn into the Jacobi-Trudi determinants $h_{\ell+1},\ldots,h_{\ell+k}$.
In description (iii), $s_\lambda$ represents the cohomology class dual
to a Schubert variety;  here $e_i$ corresponds to $s_{1^i}$, while
$h_j$ corresponds to $s_{j}$.  More generally, $s_\lambda$ can
be expressed via another Jacobi-Trudi determinant in terms of the $e_i$'s
(or in terms of the $h_j$'s).  

In what follows, we will often abuse notation by not distinguishing the indeterminate
$e_i$ in the polynomial ring $\QQ[e_1,\ldots,e_k]$ from its image in
the quotient ring $R^{k,\ell}$.

We further recall the classical expression for the {\it Hilbert series} 
$$
\begin{aligned}
\Hilb( R^{k,\ell}, q ) &:= \sum_{d \geq 0} \dim_\QQ R^{k,\ell}_d q^d \\
&= \qbinom{k+\ell}{k} 
= \frac{(1-q)(1-q^2)\cdots(1-q^{k+\ell})}
{(1-q)(1-q^2)\cdots(1-q^{k}) \cdot(1-q)(1-q^2)\cdots(1-q^{\ell})}.
\end{aligned}
$$

The first of our conjectures gives the Hilbert series for a certain
natural filtration of $R^{k,l}$.  For each $m=0,1,\ldots,k$,
let $R^{k,\ell,m}$ be the subalgebra of $R^{k,\ell}$
generated by $e_1,e_2,\ldots,e_m$, which can described alternatively
as the subalgebra generated by all elements of degree at most $m$.
Thus
$$
\QQ = R^{k,\ell,0} \subset R^{k,\ell,1} 
   \subset \cdots \subset R^{k,\ell,k} = R^{k,\ell}.
$$

\begin{conjecture}
\label{filtration-hilb}
For $m \geq 0$,
$$
\Hilb( R^{k,\ell,m}, q ) = 
  1 + \sum_{i=1}^m q^i\qbinom{\ell}{i} 
       \left( \sum_{j=0}^{k-i} q^{j(\ell-i+1)} \qbinom{i+j-1}{j}  \right).
$$
\end{conjecture}

\begin{conjecture}
\label{filtration-saturates}
For $1 \leq m \leq k$,
$$
R^{k,\ell}_d = R^{k,\ell,m-1}_d
\text{ for }d \geq k \ell - m^2 + m + 1.
$$
\end{conjecture}

\begin{conjecture}
\label{exact-e-killer}
For $1 \leq m \leq k$
$$
e_m e_1^{k \ell - m^2 + 1} \in R^{k,\ell,m-1}.
$$
\end{conjecture}

\begin{conjecture}
\label{e-overkiller}
For $3 \leq m \leq k$
$$
e_m e_1^{k \ell - 2m} \in R^{k,\ell,m-1}.
$$
\end{conjecture}

We introduce some notation for the sake of stating a result that motivated
this last conjecture.
Given a graded $\QQ$-algebra $R = \oplus_i R_i$ and $\alpha \in \QQ$,
define a graded endomorphism $\phi_\alpha$ by $\phi_\alpha(x) = \alpha^{\deg(x)} x$
for every homogeneous element $x \in R$.  In \cite{GH}, the endomorphism
$\phi_\alpha$ on $R^{k,\ell}$ for $\alpha \in \ZZ$ is called an {\it Adams map}.

Note that when $k = \ell$, there is a non-trivial automorphism
$\omega: R^{k,\ell} \rightarrow R^{k,\ell}$ (induced from the
fundamental involution on symmetric functions), exchanging
$e_i \leftrightarrow h_i$ in description (i) above, and
more generally exchanging $s_\lambda \leftrightarrow s_{\lambda'}$
in description (iii).  Here $\lambda'$ denotes the partition {\it conjugate}
to $\lambda$.

\begin{thm}(Hoffman \cite{Hoffman})
\label{non-zero-alpha}
For $k \neq \ell$, every graded algebra endomorphism 
$\phi: R^{k,\ell} \rightarrow R^{k,\ell}$ which does not annihilate
$R^{k,\ell}_1$ is of the form $\phi_\alpha$ for some $\alpha \in \QQ^\times$.

For $k = \ell$, any such endomorphism is either of the form
$\phi_\alpha$ or $\omega \circ \phi_\alpha$.
\end{thm}

Theorem~\ref{non-zero-alpha} was conjectured 
by O'Neill (see \cite{O'Neill} or \cite{GH}),
who conjectured it more generally without the assumption that 
$\phi$ is non-zero on $R^{k,\ell}_1$.  O'Neill's more general conjecture
is trivial for $k=1$ (or $\ell=1$), 
proven for $k=2$ in \cite{O'Neill}, and proven for $\ell \geq 2k^2-k-1$ in
\cite{GH}.  It is motivated by the fact that, assuming it,
one can fairly easily deduce, via the Lefschetz fixed point theorem,
that the Grassmannian $\Gr(k,\CC^{k+\ell})$ has the fixed point
property (i.e. every continuous self-map has a fixed point) if and only if $k\ell$ is
odd.

In Section~\ref{easy-implications}, we quickly explain the easy implications 
$$
\text{Conj }~\ref{filtration-hilb}
 \Rightarrow \text{Conj }~\ref{filtration-saturates} 
 \Rightarrow \text{Conj }~\ref{exact-e-killer}
 \Rightarrow \text{Conj }~\ref{e-overkiller}.
$$
In Section~\ref{easy-boundary-cases},
we verify Conjecture~\ref{filtration-hilb} in the relatively easy boundary cases
$m=1$ and $m=k$.

In Sections \ref{non-trivial-implication} and \ref{m=2-section},
we explain how Conjecture~\ref{e-overkiller} would imply
Theorem~\ref{non-zero-alpha}.  This proof uses some of the same ideas as
Hoffman's, namely the Hard Lefschetz Theorem and the hook-length formula
for counting tableaux, but is much shorter (2 pages versus 10 pages).

\section{The implications $
\text{Conj }~\ref{filtration-hilb}
 \Rightarrow \text{Conj }~\ref{filtration-saturates} 
 \Rightarrow \text{Conj }~\ref{exact-e-killer}
 \Rightarrow \text{Conj }~\ref{e-overkiller}$
}
\label{easy-implications}

To see that Conjecture~\ref{filtration-hilb} implies 
Conjecture~\ref{filtration-saturates},
note that Conjecture~\ref{filtration-hilb} is equivalent to the following
assertion: for $p \geq 1$, the quotient (graded) 
vector space $R^{k,\ell,p}/R^{k,\ell,p-1}$ has
Hilbert series
\begin{equation}
\label{filtration-quotient-hilb}
\Hilb(R^{k,\ell,p}/R^{k,\ell,p-1}, q) = \sum_{j=0}^{k-p}
q^{j(\ell-p+1)+p} \qbinom{p+j-1}{j} \cdot q^p\qbinom{\ell}{p}.
\end{equation}
It suffices to show that the right-hand has degree in $q$ at
most $k\ell - m^2 + m$ whenever $p \geq m (\geq 1)$.
Since the $q$-binomial coefficient $\qbinom{r+s}{r}$ has
degree $rs$ as a polynomial in $q$, the right-hand side of 
\eqref{filtration-quotient-hilb} has degree in $q$ equal to
$$
\begin{aligned}
\max_{j=0,1,\ldots,k-p} \{j(\ell-p+1) + p + (p-1)j + p(\ell-p)\} &=
\max_{j=0,1,\ldots,k-p} \{\ell(j+p)-p^2+p)\} \\
&= k\ell - p^2 + p,
\end{aligned}
$$
and this is bounded above by $k\ell-m^2+m$ for $p \geq m \geq 1$.

Conjecture~\ref{filtration-saturates} implies Conjecture~\ref{exact-e-killer}
trivially, since $e_m e_1^{k\ell-m^2+1}$ lies in $R^{k,\ell}_{k\ell-m^2+m+1}$.
Similarly, Conjecture~\ref{exact-e-killer}
trivially implies Conjecture~\ref{e-overkiller}, since
for $m \geq 3$, one has 
$$
k\ell-m^2+1 \leq k \ell -2m.
$$

\section{Boundary cases for Conjecture~\ref{filtration-hilb}}
\label{easy-boundary-cases}

Both in checking the boundary case $m=1$ of Conjecture~\ref{filtration-hilb}
and in later showing that Conjecture~\ref{e-overkiller} implies
Conjecture~\ref{non-zero-alpha}, we will make use of the fact that $R^{k,\ell}$
satisfies the {\it Hard Lefschetz Theorem} \cite[p. 122]{GriffithsHarris}:

\begin{thm}
\label{Hard-Lefschetz-thm}
For $i=0,1,\ldots, \lfloor \frac{k\ell}{2} \rfloor$, the map
$$
R^{k,l}_i \rightarrow R^{k,l}_{k\ell-i}
$$
given by multiplication by $e_1^{k\ell-2i}$ is a $\QQ$-vector
space isomorphism.
\end{thm}

To check the case $m=1$ in Conjecture~\ref{filtration-hilb},
note that taking $i=0$ of Theorem~\ref{Hard-Lefschetz-thm} tells us that
the smallest power $e_1^d$ which vanishes in $R^{k,\ell}$ is $e_1^{k\ell+1}$,
i.e. that $R^{k,\ell,1} \cong \QQ[e_1]/(e_1^{k\ell+1})$.
Hence 
$$
\Hilb(R^{k,\ell,1}, q) = \Hilb(\QQ[e_1]/(e_1^{k\ell+1}), q) 
= \frac{1-q^{k\ell+1}}{1-q}.
$$
Meanwhile, Conjecture~\ref{filtration-hilb} for $m=1$ predicts
$$
\Hilb(R^{k,\ell,1}, q) = 1+\sum_{j=0}^{k-1} q^{j\ell+1} \qbinom{\ell}{1}
= \frac{1-q^{k\ell+1}}{1-q},
$$
so the two agree.

\begin{remark} \rm \ \\
Instead of the Hard Lefschetz Theorem here, we could have used
Pieri's formula \cite[\S 2.2]{F} to conclude that
$$
e_1^{k\ell} = f_{\boxx} s_\boxx \neq 0
$$
where $f_\lambda$ denotes the number of {\it standard Young tableaux} of
shape $\lambda$.
\end{remark}

The case $m=k$ of Conjecture~\ref{filtration-hilb} is clearly
equivalent to the following identity.

\begin{figure}
\label{q-identity-figure}
\includegraphics{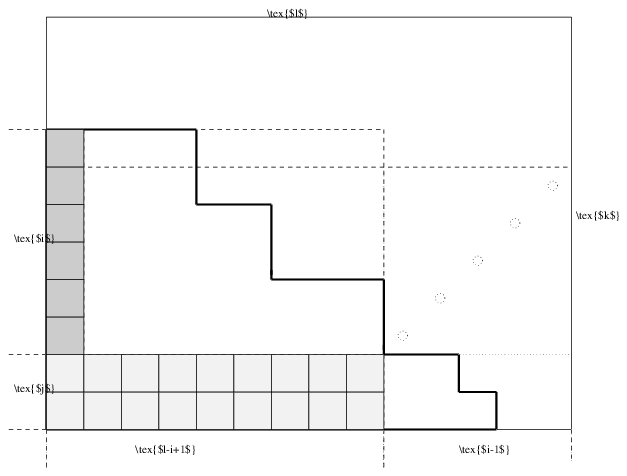}
\begin{caption}{The decomposition from the proof of Proposition~\ref{q-identity}.}
\end{caption}
\end{figure}

\begin{prop}
\label{q-identity}
$$
\qbinom{k+\ell}{k}
= 1 + 
\sum_{i=1}^k \sum_{j=0}^{k-i} q^i \qbinom{\ell}{i} 
                               \cdot q^{j(\ell-i+1)} \qbinom{i+j-1}{j} .
$$
\end{prop}
\begin{proof}
Interpret $\qbinom{k+\ell}{k}$ as the generating function counting
partitions $\lambda$ whose Ferrers diagram fits inside the rectangle $\boxx$
according to their weight.  To prove the above identity, whenever $\lambda$
is non-empty, we will uniquely define two integers $i,j$ and
decompose its diagram into four portions:
\begin{enumerate}
\item[(a)] a $j \times (\ell-i+1)$ rectangle, accounting for the
$q^{j(\ell-i+1)}$,
\item[(b)] a column of length $i$, accounting for the
$q^i$,
\item[(c)] a partition whose Ferrers diagram fits insides
a $j \times (i-1)$ rectangle, and
\item[(d)] a partition whose Ferrers diagram fits insides
a $\ell \times (\ell-i)$ rectangle.
\end{enumerate}
The decomposition is illustrated in Figure~\ref{q-identity-figure},
where Ferrers diagrams are depicted using the French notation.
Given a non-empty $\lambda$, say with $s$ parts, let 
$\bar{\lambda}$ be the partition complementary to 
$(\lambda_1,\ldots,\lambda_{s-1})$ inside a $(s-1) \times \ell$ rectangle,
and let $i-1$ be the size of the {\it Durfee square} of $\bar{\lambda}$,
that is, the largest square Ferrers diagram
contained within that of $\bar{\lambda}$.
Then set $j = s-i$.  The $j \times (\ell-i+1)$ rectangle in (a) is the one inside
$\lambda$ in the lower left.  The column of length $i$ in (b) lies just
above it in column $1$.  The remaining cells of $\lambda$ then
segregate into two Ferrers diagrams, one inside a $j \times (i-1)$ rectangle
(as in (c)) in the lower right, the other inside an $\ell \times (\ell-i)$ rectangle
(as in (d)) in the upper left.
\end{proof}

\begin{remark} \rm \ \\
It may be that a more useful phrasing for Conjecture~\ref{filtration-hilb}
is to rewrite the inner sum in the following form, which one can show is equivalent:
$$
\Hilb( R^{k,\ell,m}, q ) = 
  1 + \sum_{i=1}^m q^i\qbinom{\ell}{i} f^{k,\ell}_i(q)
$$
where $f^{k,\ell}_i(q)$ for $0 \leq i \leq k \leq \ell$
is defined by the following recurrence
$$
f^{k,\ell}_i(q) = f^{k-1,\ell}_{i-1}(q) + q^{\ell-i+1} f^{k-1,\ell}_i(q)
$$
with initial conditions $f^{k,\ell}_0(q) = f^{k,\ell}_k(q)=1$.
In other words, $f^{k,\ell}_i(q)$ is a $q$-analogue of the binomial
coefficient $\binom{k}{i}$ which depends on $\ell$ also, and satisfies
a different $q$-Pascal recurrence than the usual one for the
$q$-analogue $\qbinom{k}{i}$.
\end{remark}

\section{ Conjecture~\ref{e-overkiller} implies Theorem~\ref{non-zero-alpha}
}
\label{non-trivial-implication}

Let $\phi: R^{k,\ell} \rightarrow R^{k,\ell}$ be a graded algebra
endomorphism which does not annihilate $R^{k,\ell}_1$.  Since
$R^{k,\ell}_1$ is a $1$-dimensional $\QQ$-vector space, 
$\phi$ acts on $R^{k,\ell}_1$ by some constant $\alpha \in \QQ^\times$.
Recall that $\phi_\alpha$ is the endomorphism that sends a homogeneous
element $x$ in $R^{k,\ell}_d$ to $\alpha^d x$. 
It is our goal in this (and the next) 
section to show that (assuming Conjecture~\ref{e-overkiller})
\begin{enumerate}
\item[$\bullet$]
if $k \neq \ell$ then
$\phi = \phi_\alpha$, and 
\item[$\bullet$]
if $k = \ell$ either 
$\phi = \phi_\alpha$ or $\phi = \omega \circ \phi_\alpha$.
\end{enumerate}

Note that $\phi_{\alpha^{-1}} = \phi_\alpha^{-1}$, and so
by replacing $\phi$ with the composite $\phi \circ \phi_{\alpha^{-1}}$,
we may assume without loss of generality that $\alpha=1$, i.e.
$\phi$ acts as the identity on $R^{k,\ell}_1$.

Our goal will be to try and show by induction on $m$ that
$\phi$ does what we expect to $e_1,e_2,\ldots,m$ for every $m \leq k$.
The case $m=2$ is summarized in the following lemma, whose
proof requires ad hoc argumenation which we defer to Section~\ref{m=2-section}.

\begin{lemma}
\label{m=2}
Let $\phi: R^{k,\ell} \rightarrow R^{k,\ell}$ 
be a graded algebra endomorphism which acts as the identity
on $R^{k,\ell}_1$.

If $k \neq \ell$, then $\phi(e_2) = e_2$, i.e. $\phi$ acts the
identity on $R^{k,\ell}_2$.

If $k = \ell$, then $\phi(e_2)$ is either $e_2$ or $h_2$, i.e.
$\phi$ acts either as the identity or $\omega$ on $R^{k,\ell}_2$.
\end{lemma}

Assuming this lemma for the moment, we show how 
Conjecture~\ref{e-overkiller} allows one to do the remaining inductive
steps for $m \geq 3$ to deduce Conjecture~\ref{non-zero-alpha}.

Given $m \geq 3$, using Lemma~\ref{m=2}, we may assume without
loss of generality (by composing $\phi$ with $\omega$ if necessary when $k = \ell$),
that $\phi(e_i)=e_i$ for $i=1,2,\ldots,m-1$
by induction on $m$.  In other words, $\phi$ acts as the identity on the
subalgebra $R^{k,\ell,m-1}$.
We wish to show that this implies $\phi(e_m)=e_m$.

Let $\PP_{k,\ell}(m)$ denote the set of all partitions $\lambda$ of $m$
whose Ferrers diagram fits inside the rectangle $\boxx$.  Note that
for all but one such $\lambda$, namely $\lambda=(m)$ having a single
part, the elementary symmetric function 
$e_\lambda: = e_{\lambda_1} e_{\lambda_2} \cdots e_{\lambda_k}$ 
lies in the subalgebra $R^{k,\ell,m-1}$, and hence so
does the product $e_\lambda e_1^{k\ell-2m}$.  
{\it Because of Conjecture~\ref{e-overkiller}}, the same is true for
$\lambda = 1^m$, i.e. 
$$
e_{(m)} e_1^{k\ell-2m} =  e_m e_1^{k\ell-2m} \in R^{k,\ell,m-1}.
$$
This is the only place where we will use our assumption of Conjecture~\ref{e-overkiller}.

Hence $\phi$ fixes $e_\lambda e_1^{k\ell-2m}$ for every 
$\lambda \in \PP_{k,\ell}(m)$.  From degree considerations, there exists
constants $b_\lambda \in \QQ$ for $\lambda \in \PP_{k,\ell}(m)$
(whose exact values are not important) such that
\begin{equation}
\label{e-multiplication}
e_m e_\lambda e_1^{k\ell-2m} = b_\lambda s_\boxx.
\end{equation}
Let 
$$
\phi(e_m) = \sum_{\mu \in \PP_{k,\ell}(m)} x_\mu s_\mu
$$
for some constants $x_\mu \in \QQ$.
We hope to show that $x_\mu = \delta_{\mu, (m)}$, so that $\phi$ fixes $e_m$.
Applying $\phi$
to \eqref{e-multiplication} (and noting that $\phi$ also fixes
$s_\boxx$ because it lies in $R^{k,\ell,1}$) yields the system of equations
\begin{equation}
\sum_{\mu \in \PP_{k,\ell}(m)} x_\mu s_\mu e_\lambda e_1^{k\ell-2m} = 
b_\lambda s_\boxx \text{ for }\lambda \in \PP_{k,\ell}(m).
\end{equation}
This system can be written in matrix form as $Ax=b$, where
$A=(a_{\lambda,\mu})_{\lambda,\mu \in \PP_{k,\ell}(m)}$,
and $a_{\lambda,\mu}$ is the coefficient of $s_\boxx$ when
one expands $s_\mu e_\lambda e_1^{k\ell-2m}$.

Now we know that this system has at least one solution 
$x_\mu = \delta_{\mu, (m)}$, corresponding to the identity
endomorphism $\phi$.  To show this is the unique solution,
we need only show that $A$ is invertible.  But this follows from the
Hard Lefschetz Theorem~\ref{Hard-Lefschetz-thm} and Poincar\'e
duality, as the matrix $A$ expresses the (invertible) linear map 
$$
R^{k,\ell}_m \qquad 
\overset{\cdot e_1^{k\ell-2m}}{\longrightarrow} \qquad R^{k,\ell}_{k\ell-m}
$$
with respect to a particular choice of $\QQ$-bases in the domain, range:
\begin{enumerate}
\item[$\bullet$]
In the domain, the basis $\{ e_\lambda \}_{\lambda \in \PP_{k,\ell}(m)}$ for
$R^{k,\ell}_m$.  This is a basis because it is upper-triangularly related to the 
usual basis $\{ s_\lambda \}_{\lambda \in \PP_{k,\ell}(m)}$.
\item[$\bullet$]
In the range, the basis for $R^{k,\ell}_{k\ell-m}$ which is {\it Poincar\'e dual} to 
the basis $\{ s_\mu \}_{\mu \in \PP_{k,\ell}(m)}$ for $R^{k,\ell}_m$.
\end{enumerate}

This completes the inductive step when $m \geq 3$, and leaves
only the proof of Lemma~\ref{m=2} remaining.

\begin{remark} \rm \ \\
We note a further consequence of the Hard Lefschetz Theorem~\ref{Hard-Lefschetz-thm}
here.  It implies the coincidence of $\QQ$-vector spaces  
$$
e_1^{k\ell-2m} R^{k,\ell}_m = R^{k,\ell}_{k\ell-m}.
$$
Since $e_m$ is the only element of the basis 
$\{ e_\lambda \}_{\lambda \in \PP_{k,\ell}(m)}$ for
$R^{k,\ell}_m$ which does not lie in the subalgebra $R^{k,\ell,m-1}$,
from this one sees that Conjecture~\ref{e-overkiller} is
equivalent to the following statement, a weakening
of Conjecture~\ref{filtration-saturates}:

\vskip .1in
\noindent
{\bf Conjecture ~\ref{e-overkiller}$^\prime$.}
{\it
For $3 \leq m \leq k$,
$$
R^{k,\ell}_d = R^{k,\ell,m-1}_d
\text{ for }d \geq k \ell - 2m.
$$
}
\end{remark}

\section{Proof of Lemma~\ref{m=2}.}
\label{m=2-section}

We recall the statement of the lemma.

\noindent
{\bf Lemma~\ref{m=2}.}
{\it 
Let $\phi: R^{k,\ell} \rightarrow R^{k,\ell}$ 
be a graded algebra endomorphism which acts as the identity
on $R^{k,\ell}_1$.

If $k \neq \ell$, then $\phi(e_2) = e_2$, i.e. $\phi$ acts the
identity on $R^{k,\ell}_2$.

If $k = \ell$, then $\phi(e_2)$ is either $e_2$ or $h_2$, i.e.
$\phi$ acts either as the identity or $\omega$ on $R^{k,\ell}_2$.
}

\begin{proof}
Assume that $\phi$ acts as the identity on $R^{k,\ell}_1$,
so $\phi(e_1)=e_1$, and let $\phi(e_2) = x e_2 + y e_1^2$
for two unknown constants $x, y \in \QQ$.

From degree considerations, there are constants $\gamma_r \in \QQ$ for
for $r=0,1,2,\ldots$ satisfying
\begin{equation}
\label{e2-equations}
e_2^r e_1^{k \ell - 2r} = \gamma_r s_\boxx.
\end{equation}
For small values of $r$ one can use the Pieri formula \cite[\S2.2]{F}
to obtain explicit formulae for these constants:
$$
\begin{aligned}
\gamma_0 &= f_{\emptyset^c} \\
\gamma_1 &= f_{11^c} \\
\gamma_2 &= f_{1111^c} + f_{211^c} + f_{22^c}\\
\gamma_3 &= f_{11111^c} + 2f_{21111^c} + 3f_{2211^c} + f_{222^c} + f_{3111^c} + 
  2f_{321^c} + f_{33^c}
\end{aligned}
$$
where here $\lambda^c$ denotes the partition whose Ferrers diagram is the
complement of that of $\lambda$ within the rectangle $\boxx$, and we recall
that $f_\mu$ denotes number of {\it standard Young tableaux} of shape $\mu$.
One can apply the endomorphism $\phi$ to the equations \eqref{e2-equations},
(recalling that $\phi$ fixes both $e_1$ and $s_\boxx$ because they
lie in $R^{k,\ell,1}$), expand using the expression for $\phi(e_2)$,
and then divide both sides by $\gamma_0$,
to obtain the following equations for $r = 0,1,2,\ldots$
in the unknowns $x, y$:
\begin{equation}
\label{e2-derived-equations}
\sum_{i=0}^r \binom{r}{i}\frac{\gamma_i}{\gamma_0} x^i y^{r-i} 
= \frac{\gamma_r}{\gamma_0}.
\end{equation}

Here the celebrated {\it hook-length formula} \cite[\S 4.3]{F} for $f_\lambda$
comes to the rescue: for $\lambda$ a partition of $n$, one has
$$
f_\lambda = \frac{n!}{\prod_{x \in \lambda} h(x) }
$$
where the product is taken over all cells $x$ in the Ferrers diagram for
$\lambda$, and $h(x)$ is the {\it hooklength} at cell $x$, that is, the
number of cells in the diagram that are either weakly to the right of
$x$ in the same row or weakly below it in the same column.
Using this formula, the constants $\frac{\gamma_r}{\gamma_0}$ for small
values of $r$ can be explicitly computed as rational functions of $k, \ell$,
with relatively small numerators and denominators;  we omit these formulae here.

With the aid of these formulae (and some help from computer
algebra packages), one can check that the only simultaneous solutions
to the two equations \eqref{e2-derived-equations} for $r=1,2$
are
\begin{enumerate}
\item[(i)] $(x,y) =(1,0)$, and
\item[(ii)] $(x,y) = (-1, \frac{(k-1)(\ell+1)}{k\ell-1}).$
\end{enumerate}
Solution (i) gives $\phi(e_2) = e_2$, which is a scenario that we hoped
for. Substituting solution (ii) into
\eqref{e2-derived-equations} for $r=3$ gives
$$
\frac{(\ell+1)(k-1)(k+1)(\ell-1)(k\ell+5)(k-\ell)}
         {(k\ell-1)^3(k\ell-2)(k\ell-3)(k\ell-4)(k\ell-5)}=0.
$$
Since we may assume without loss of generality that $k, \ell \geq 2$,
this forces $k=\ell$.  In this case, solution (ii) becomes
$(x,y) = (-1, 1)$, which means $\phi(e_2)=-e_1^2+e_2 = h_2$, as desired.

\end{proof}

\end{document}